\documentclass[12pt,a4paper]{article}
\usepackage{fullpage}
\usepackage{amsthm}
\usepackage{amsfonts}
\usepackage{amsmath}
\usepackage{amscd}
\usepackage{amssymb}
\usepackage{pstricks,pst-node,pst-tree}
\usepackage{stmaryrd}
\theoremstyle{plain}
\newtheorem{thm}{Theorem}[section]
\newtheorem{lem}[thm]{Lemma}
\newtheorem{prop}[thm]{Proposition}
\newtheorem{cor}[thm]{Corollary}

\newtheorem*{nathanson}{Theorem~3.6 of \cite{Nathanson2000}}

\theoremstyle{definition}

\theoremstyle{remark}
\newtheorem*{rem}{Remark}

\newcommand{\N}{\mathbb{N}}
\newcommand{\Z}{\mathbb{Z}}
\newcommand{\Zn}{\Z/n\Z}

\newcommand{\ord}{\mathcal{O}}
\newcommand{\pord}{\mathcal{PO}}
\newcommand{\Rn}{\mathcal{R}_n}
\newcommand{\twoheadlongrightarrow}{\relbar\joinrel\twoheadrightarrow}
\DeclareMathOperator*{\lcm}{lcm}
\DeclareMathOperator*{\rad}{rad}
\renewcommand{\leq}{\leqslant}
\renewcommand{\geq}{\geqslant}
\title{On the multiplicative order of $a^n$ modulo $n$}
\date{January 20, 2010}
\author{Jonathan Chappelon\textsuperscript{a,b,c}\\
\noindent
\textsuperscript{a}Univ Lille Nord de France, F-59000 Lille, France\\
\textsuperscript{b}ULCO, LMPA J.~Liouville, B.P. 699, F-62228 Calais, France\\
\textsuperscript{c}CNRS, FR 2956, France }
\begin{document}
\maketitle
\begin{abstract}
Let $n$ be a positive integer and $\alpha_n$ be the arithmetic function which assigns the multiplicative order of $a^n$ modulo $n$ to every integer $a$ coprime to $n$ and vanishes elsewhere. Similarly, let $\beta_n$ assign the projective multiplicative order of $a^n$ modulo $n$ to every integer $a$ coprime to $n$ and vanishes elsewhere. In this paper, we present a study of these two arithmetic functions. In particular, we prove that for positive integers $n_1$ and $n_2$ with the same square-free part, there exists an exact relationship between the functions $\alpha_{n_1}$ and $\alpha_{n_2}$ and between the functions $\beta_{n_1}$ and $\beta_{n_2}$. This allows us to reduce the determination of $\alpha_n$ and $\beta_n$ to the case where $n$ is square-free. These arithmetic functions recently appeared in the context of an old problem of Molluzzo, and more precisely in the study of which arithmetic progressions yield a balanced Steinhaus triangle in $\Zn$ for $n$ odd.\\[2ex]
\end{abstract}
\section{Introduction}
We start by introducing some notation relating to the order of certain elements modulo $n$. For every positive integer $n$ and every prime number $p$, we denote by $v_p(n)$ the \textit{$p$-adic valuation} of $n$, i.e., the greatest exponent $e\geq0$ for which $p^e$ divides $n$. The prime factorization of $n$ may then be written as
$$
n=\prod_{p\in\mathbb{P}}p^{v_p(n)},
$$
where $\mathbb{P}$ denotes the set of all prime numbers. We denote by $\rad(n)$ \textit{the radical of $n$}, i.e., the largest square-free divisor of $n$, namely
$$
\rad(n) = \prod_{\stackrel{p\in\mathbb{P}}{p|n}}p.
$$
\par For every positive integer $n$ and every integer $a$ coprime to $n$, we denote by $\ord_n(a)$ \textit{the multiplicative order of $a$ modulo $n$}, i.e., the smallest positive integer $e$ such that $a^e\equiv1\pmod{n}$, namely
$$
\ord_n(a) = \min\left\{e\in\N^*\ \middle|\ a^e\equiv1\pmod{n}\right\}
$$
and we denote by $\Rn(a)$ the \textit{multiplicative remainder} of $a$ modulo $n$, i.e., the multiple of $n$ defined by
$$
\Rn(a) = a^{\ord_n(a)}-1.
$$
The multiplicative order of $a$ modulo $n$ also corresponds with the order of the element $\pi_n(a)$, where $\pi_n : \Z \twoheadlongrightarrow \Zn$ is the canonical surjective morphism, in the multiplicative group ${\left(\Zn\right)}^*$, the group of units of $\Zn$. Note that $\ord_n(a)$ divides $\varphi(n)$, $\varphi$ being the Euler's totient function.
\par For every positive integer $n$, we define and denote by $\alpha_n$ the arithmetic function
$$
\alpha_n :
\begin{array}[t]{ccl}
\Z & \longrightarrow & \N\\
a & \longmapsto & \left\{
\begin{array}{cl}
\ord_n\left(a^n\right), & \text{for\ all}\ \gcd(a,n)=1;\\
0, & \text{otherwise},
\end{array}\right.
\end{array}
$$
where $\gcd(a,n)$ denotes the greatest common divisor of $a$ and $n$, with the convention that $\gcd(0,n)=n$. Observe that, for every $a$ coprime to $n$, the integer $\alpha_n(a)$ divides $\varphi(n)/\gcd(\varphi(n),n)$. This follows from the previous remark on $\ord_n(a)$ and the equality $\alpha_n(a)=\ord_n(a^n)=\ord_n(a)/\gcd(\ord_n(a),n)$.
\par For every positive integer $n$ and every integer $a$ coprime to $n$, we denote by $\pord_n(a)$ \textit{the projective multiplicative order of $a$ modulo $n$}, i.e., the smallest positive integer $e$ such that $a^e\equiv\pm1\pmod{n}$, namely
$$
\pord_n(a) = \min\left\{e\in\N^*\ \middle|\ a^e\equiv\pm1\pmod{n}\right\}.
$$
The projective multiplicative order of $a$ modulo $n$ also corresponds with the order of the element $\overline{\pi_n(a)}$ in the multiplicative quotient group ${\left(\Zn\right)}^*/\{-1,1\}$.
\par For every positive integer $n$, we define and denote by $\beta_n$ the arithmetic function
$$
\beta_n :
\begin{array}[t]{ccl}
\Z & \longrightarrow & \N\\
a & \longmapsto & \left\{
\begin{array}{cl}
\pord_n\left(a^n\right), & \text{for\ all}\ \gcd(a,n)=1;\\
0, & \text{otherwise}.
\end{array}\right.
\end{array}
$$
Observe that we have the alternative $\alpha_n=\beta_n$ or $\alpha_n=2\beta_n$.
\par
\begin{figure}[!h]
\begin{center}
\begin{pspicture}(3,1.9485)
\pspolygon(0.375,1.9485)(2.625,1.9485)(1.5,0)
\rput(1.5,0.433){$0$}
\rput(1.25,0.866){$4$}
\rput(1.75,0.866){$1$}
\rput(1,1.299){$4$}
\rput(1.5,1.299){$0$}
\rput(2,1.299){$1$}
\rput(0.75,1.732){$2$}
\rput(1.25,1.732){$2$}
\rput(1.75,1.732){$3$}
\rput(2.25,1.732){$3$}
\end{pspicture}
\end{center}
\vspace{-0.75cm}
\caption{\label{fig2}The Steinhaus triangle of a balanced sequence in $\Z/5\Z$}
\end{figure}
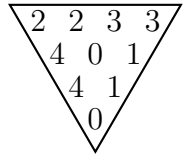

 In this paper, we study in detail these two arithmetic functions. In particular, we prove that, for every positive integers $n_1$ and $n_2$ such that
$$
\left\{
\begin{array}{ll}
\rad(n_1)|n_2\ \text{and}\ n_2|n_1, & \text{if}\ v_2(n_1)\leq1;\\
2\rad(n_1)|n_2\ \text{and}\ n_2|n_1, & \text{if}\  v_2(n_1)\geq2,
\end{array}
\right.
$$
the integer $\alpha_{n_1}(a)$ (respectively $\beta_{n_1}(a)$) divides $\alpha_{n_2}(a)$ (resp. $\beta_{n_2}(a)$), for every integer $a$. More precisely, we determine the exact relationship between the functions $\alpha_{n_1}$ and $\alpha_{n_2}$ and between $\beta_{n_1}$ and $\beta_{n_2}$. We prove that we have
$$
\alpha_{n_1}(a) = \frac{\alpha_{n_2}(a)}{\gcd\left(\alpha_{n_2}(a),\frac{\gcd(n_1,\mathcal{R}_{n_2}(a))}{n_2}\right)}\ \text{for\ all}\ \gcd(a,n_1)=1
$$
in Theorem~\ref{thm21} of Section 2 and that we have
$$
\beta_{n_1}(a) = \frac{\beta_{n_2}(a)}{\gcd\left(\beta_{n_2}(a),\frac{\gcd(n_1,\mathcal{R}_{n_2}(a))}{n_2}\right)}\ \text{for\ all}\ \gcd(a,n_1)=1
$$
in Theorem~\ref{thm31} of Section 3. Thus, for every integer $a$ coprime to $n$, the determination of $\alpha_n(a)$ is reduced to the computation of $\alpha_{\rad(n)}(a)$ and $\mathcal{R}_{\rad(n)}(a)$ if $v_2(n)\leq1$ and of $\alpha_{2\rad(n)}(a)$ and $\mathcal{R}_{2\rad(n)}(a)$ if $v_2(n)\geq2$. These theorems on the functions $\alpha_n$ and $\beta_n$ are derived from Theorem~\ref{thm2} of Section~2, which states that
$$
\ord_{n_1}(a) = \ord_{n_2}(a)\cdot\frac{n_1}{\gcd(n_1,\mathcal{R}_{n_2}(a))},
$$
for all integers $a$ coprime to $n_1$ and $n_2$. This result generalizes the following theorem of Nathanson which, in the above notation, states that for every odd prime number $p$ and for every positive integer $k$, we have the equality
$$
\ord_{p^k}(a)=\ord_p(a)\cdot\frac{p^k}{\gcd(p^k,\mathcal{R}_p(a))}
$$
for all integers $a$ not divisible by $p$.
\begin{nathanson}
Let $p$ be an odd prime, and let $a\neq\pm1$ be an integer not divisible by $p$. Let $d$ be the order of $a$ modulo $p$. Let $k_0$ be the largest integer such that $a^d\equiv1\pmod{p^{k_0}}$. Then the order of $a$ modulo $p^k$ is $d$ for $k=1,\ldots,k_0$ and $dp^{k-k_0}$ for $k\geq k_0$.
\end{nathanson}
For every finite sequence $S=\left(a_1,\ldots,a_m\right)$ of length $m\geq1$ in $\Zn$, we denote by $\Delta S$ the \textit{Steinhaus triangle} of $S$, that is the finite multiset of cardinality $\binom{m+1}{2}$ in $\Zn$ defined by
$$
\Delta S = \left\{ \sum_{k=0}^{i}{\binom{i}{k}a_{j+k}}\ \middle|\ 0\leq i\leq m-1\ ,\ 1\leq j\leq m-i \right\}.
$$
A finite sequence $S$ in $\Zn$ is said to be \textit{balanced} if each element of $\Zn$ occurs in its Steinhaus triangle $\Delta S$ with the same multiplicity. For instance, the sequence $\left(2,2,3,3\right)$ of length $4$ is balanced in $\Z/5\Z$. Indeed, as depicted in Figure~\ref{fig2}, its Steinhaus triangle is composed by each element of $\Z/5\Z$ occuring twice.

Note that, for a sequence $S$ of length $m\geq1$ in $\Zn$, a necessary condition on $m$ for $S$ to be balanced is that the integer $n$ divides the binomial coefficient $\binom{m+1}{2}$. In 1976, John C. Molluzzo \cite{Molluzzo1978} posed the problem to determine whether this necessary condition on $m$ is also sufficient to guarantee the existence of a balanced sequence. In \cite{Chappelon}, it was proved that, for each odd number $n$, \textit{there exists a balanced sequence of length $m$ for every $m\equiv0$ or $-1\pmod{\alpha_n(2)\cdot n}$ and for every $m\equiv0$ or $-1\pmod{\beta_n(2)\cdot n}$}. This was achieved by analyzing the Steinhaus triangles generated by arithmetic progressions. In particular, since $\beta_{3^k}(2)=1$ for all $k\geq1$, the above result implies a complete and positive solution of Molluzzo's Problem in $\Zn$ for all $n=3^k$.

\section{The arithmetic function $\alpha_n$}

The table depicted in Figure~\ref{fig3} gives us the first values of $\alpha_n(a)$ for every positive integer $n$, $1\leq n\leq 20$, and for every integer $a$, $-20\leq a\leq 20$.

\begin{figure}[!h]
\begin{center}
\begin{tabular}{|c||c|c|c|c|c|c|c|c|c|c|c|c|c|c|c|c|c|c|c|c|}
\hline
$n\backslash a$ & $1$ & $2$ & $3$ & $4$ & $5$ & $6$ & $7$ & $8$ & $9$ & $10$ & $11$ & $12$ & $13$ & $14$ & $15$ & $16$ & $17$ & $18$ & $19$ & $20$\\
\hline
\hline
$1$ & $1$ & $1$ & $1$ & $1$ & $1$ & $1$ & $1$ & $1$ & $1$ & $1$ & $1$ & $1$ & $1$ & $1$ & $1$ & $1$ & $1$ & $1$ & $1$ & $1$\\
\hline
$2$ & $1$ & $0$ & $1$ & $0$ & $1$ & $0$ & $1$ & $0$ & $1$ & $0$ & $1$ & $0$ & $1$ & $0$ & $1$ & $0$ & $1$ & $0$ & $1$ & $0$\\
\hline
$3$ & $1$ & $2$ & $0$ & $1$ & $2$ & $0$ & $1$ & $2$ & $0$ & $1$ & $2$ & $0$ & $1$ & $2$ & $0$ & $1$ & $2$ & $0$ & $1$ & $2$\\
\hline
$4$ & $1$ & $0$ & $1$ & $0$ & $1$ & $0$ & $1$ & $0$ & $1$ & $0$ & $1$ & $0$ & $1$ & $0$ & $1$ & $0$ & $1$ & $0$ & $1$ & $0$\\
\hline
$5$ & $1$ & $4$ & $4$ & $2$ & $0$ & $1$ & $4$ & $4$ & $2$ & $0$ & $1$ & $4$ & $4$ & $2$ & $0$ & $1$ & $4$ & $4$ & $2$ & $0$\\
\hline
$6$ & $1$ & $0$ & $0$ & $0$ & $1$ & $0$ & $1$ & $0$ & $0$ & $0$ & $1$ & $0$ & $1$ & $0$ & $0$ & $0$ & $1$ & $0$ & $1$ & $0$\\
\hline
$7$ & $1$ & $3$ & $6$ & $3$ & $6$ & $2$ & $0$ & $1$ & $3$ & $6$ & $3$ & $6$ & $2$ & $0$ & $1$ & $3$ & $6$ & $3$ & $6$ & $2$\\
\hline
$8$ & $1$ & $0$ & $1$ & $0$ & $1$ & $0$ & $1$ & $0$ & $1$ & $0$ & $1$ & $0$ & $1$ & $0$ & $1$ & $0$ & $1$ & $0$ & $1$ & $0$\\
\hline
$9$ & $1$ & $2$ & $0$ & $1$ & $2$ & $0$ & $1$ & $2$ & $0$ & $1$ & $2$ & $0$ & $1$ & $2$ & $0$ & $1$ & $2$ & $0$ & $1$ & $2$\\
\hline
$10$ & $1$ & $0$ & $2$ & $0$ & $0$ & $0$ & $2$ & $0$ & $1$ & $0$ & $1$ & $0$ & $2$ & $0$ & $0$ & $0$ & $2$ & $0$ & $1$ & $0$\\
\hline
$11$ & $1$ & $10$ & $5$ & $5$ & $5$ & $10$ & $10$ & $10$ & $5$ & $2$ & $0$ & $1$ & $10$ & $5$ & $5$ & $5$ & $10$ & $10$ & $10$ & $5$\\
\hline
$12$ & $1$ & $0$ & $0$ & $0$ & $1$ & $0$ & $1$ & $0$ & $0$ & $0$ & $1$ & $0$ & $1$ & $0$ & $0$ & $0$ & $1$ & $0$ & $1$ & $0$\\
\hline
$13$ & $1$ & $12$ & $3$ & $6$ & $4$ & $12$ & $12$ & $4$ & $3$ & $6$ & $12$ & $2$ & $0$ & $1$ & $12$ & $3$ & $6$ & $4$ & $12$ & $12$\\
\hline
$14$ & $1$ & $0$ & $3$ & $0$ & $3$ & $0$ & $0$ & $0$ & $3$ & $0$ & $3$ & $0$ & $1$ & $0$ & $1$ & $0$ & $3$ & $0$ & $3$ & $0$\\
\hline
$15$ & $1$ & $4$ & $0$ & $2$ & $0$ & $0$ & $4$ & $4$ & $0$ & $0$ & $2$ & $0$ & $4$ & $2$ & $0$ & $1$ & $4$ & $0$ & $2$ & $0$\\
\hline
$16$ & $1$ & $0$ & $1$ & $0$ & $1$ & $0$ & $1$ & $0$ & $1$ & $0$ & $1$ & $0$ & $1$ & $0$ & $1$ & $0$ & $1$ & $0$ & $1$ & $0$\\
\hline
$17$ & $1$ & $8$ & $16$ & $4$ & $16$ & $16$ & $16$ & $8$ & $8$ & $16$ & $16$ & $16$ & $4$ & $16$ & $8$ & $2$ & $0$ & $1$ & $8$ & $16$\\
\hline
$18$ & $1$ & $0$ & $0$ & $0$ & $1$ & $0$ & $1$ & $0$ & $0$ & $0$ & $1$ & $0$ & $1$ & $0$ & $0$ & $0$ & $1$ & $0$ & $1$ & $0$\\
\hline
$19$ & $1$ & $18$ & $18$ & $9$ & $9$ & $9$ & $3$ & $6$ & $9$ & $18$ & $3$ & $6$ & $18$ & $18$ & $18$ & $9$ & $9$ & $2$ & $0$ & $1$\\
\hline
$20$ & $1$ & $0$ & $1$ & $0$ & $0$ & $0$ & $1$ & $0$ & $1$ & $0$ & $1$ & $0$ & $1$ & $0$ & $0$ & $0$ & $1$ & $0$ & $1$ & $0$\\
\hline
\end{tabular}
\end{center}
\caption{\label{fig3}The first values of $\alpha_n(a)$}
\end{figure}

The positive integer $\alpha_n(a)$ seems to be difficult to determine. Indeed, there is no general formula known to compute the multiplicative order of an integer modulo $n$ but, however, we get the following helpful propositions.

\begin{lem}\label{lem11}
Let $n_1$ and $n_2$ be two positive integers such that $\rad(n_1)=\rad(n_2)$. Then, an integer $a$ is coprime to $n_1$ if, and only if, it is also coprime to $n_2$.
\end{lem}

\begin{proof}
This follows from the definition of the greatest common divisor of two integers and from the definition of the radical of an integer.
\end{proof}

\begin{prop}\label{prop2}
Let $n_1$ and $n_2$ be two positive integers such that $\rad(n_1)|n_2$ and $n_2|n_1$. Then, for every integer $a$, the integer $\alpha_{n_1}(a)$ divides $\alpha_{n_2}(a)$.
\end{prop}

\begin{proof}
If $a$ is not coprime to $n_1$ and $n_2$, then, by definition of the functions $\alpha_{n_1}$ and $\alpha_{n_2}$ and by Lemma~\ref{lem11}, we have
$$
\alpha_{n_1}(a)=\alpha_{n_2}(a)=0.
$$
Suppose that $a$ is coprime to $n_1$ and $n_2$. If $v_p(n_1)=1$ for all prime factors $p$ of $n_1$, then $n_2=n_1$. Otherwise, let $p$ be a prime factor of $n_1$ such that $v_p(n_1)\geq2$. We shall show that $\alpha_{n_1}(a)$ divides $\alpha_{n_1/p}(a)$. By definition of $\alpha_{n_1/p}(a)$, there exists an integer $u$ such that
$$
a^{\alpha_{n_1/p}(a)\cdot\frac{n_1}{p}} = 1+u\cdot\frac{n_1}{p}.
$$
Therefore, by the binomial theorem, we have
$$
a^{\alpha_{n_1/p}(a)\cdot n_1} = {\left(a^{\alpha_{n_1/p}(a)\cdot\frac{n_1}{p}}\right)}^{p} = {\left(1+u\cdot\frac{n_1}{p}\right)}^{p} = 1+ u\cdot n_1 + \sum_{k=2}^{p}{\binom{p}{k}\cdot u^k\cdot\left(\frac{n_1}{p}\right)^k}.
$$
Since $v_p(n_1)\geq2$, it follows that ${\left(n_1/p\right)}^k$ is divisible by $n_1$ for every integer $k\geq2$ and so
$$
a^{\alpha_{n_1/p}(a)\cdot n_1} \equiv 1 \pmod{n_1}.
$$
Hence $\alpha_{n_1}(a)$ divides $\alpha_{n_1/p}(a)$. This completes the proof.
\end{proof}

An exact relationship between $\alpha_{n_1}(a)$ and $\alpha_{n_2}(a)$, for every integer $a$ coprime to $n_1$ and $n_2$, is determined at the end of this section. We first settle the easy prime power case.

\begin{prop}\label{prop3}
Let $p$ be a prime number and let $a$ be an integer. Then we have
$$
\alpha_{p^k}(a) = \ord_p(a)
$$
for every positive integer $k$.
\end{prop}

\begin{proof}
Let $k$ be a positive integer. If $a$ is not coprime to $p$, then  we have $\alpha_{p^k}(a)=\alpha_{p}(a)=0$. Suppose now that $a$ is coprime to $p$. By Proposition~\ref{prop2}, the integer $\alpha_{p^k}(a)$ divides $\alpha_{p}(a)$. It remains to prove that $\alpha_{p}(a)$ divides $\alpha_{p^k}(a)$. The congruence
$$
a^{\alpha_{p^k}(a)\cdot p^k} \equiv 1 \pmod{p^k}
$$
implies that
$$
a^{\alpha_{p^k}(a)\cdot p^k} \equiv 1 \pmod{p},
$$
and hence, by Fermat's Little Theorem, it follows that
$$
a^{\alpha_{p^k}(a)\cdot p} \equiv a^{\alpha_{p^k}(a)\cdot p^k} \equiv 1 \pmod{p}.
$$
Therefore $\alpha_{p}(a)$ divides $\alpha_{p^k}(a)$. Finally, we have
$$
\alpha_{p^k}(a) = \alpha_{p}(a) = \ord_{p}(a^p) = \ord_{p}(a).
$$
This completes the proof.
\end{proof}

\begin{rem}
If $p=2$, then, for every positive integer $k$, we obtain
$$
\alpha_{2^k}(a) = \ord_{2}(a) = \left\{
\begin{array}{ll}
0, & \text{for}\ a\ \text{even};\\
1, & \text{for}\ a\ \text{odd}.
\end{array}
\right.
$$
\end{rem}

\begin{prop}
Let $n_1$ and $n_2$ be two coprime numbers and let $a$ be an integer. Then $\alpha_{n_1n_2}(a)$ divides $\lcm(\alpha_{n_1}(a),\alpha_{n_2}(a))$, the least common multiple of $\alpha_{n_1}(a)$ and $\alpha_{n_2}(a)$.
\end{prop}

\begin{proof}
If $\gcd(a,n_1n_2)\neq1$, then $\gcd(a,n_1)\neq1$ or $\gcd(a,n_2)\neq1$ and so
$$
\alpha_{n_1n_2}(a) = \lcm(\alpha_{n_1}(a),\alpha_{n_2}(a)) = 0.
$$
Suppose now that $\gcd(a,n_1n_2)=1$ and hence that the integers $a$, $n_1$ and $n_2$ are coprime pairwise. Let $i\in\{1,2\}$. The congruences
$$
a^{\alpha_{n_i}(a)\cdot n_i}\equiv1\pmod{n_i}
$$
imply that
$$
a^{n_1n_2\lcm(\alpha_{n_1}(a),\alpha_{n_2}(a))}\equiv1\pmod{n_i}.
$$
Therefore $\alpha_{n_1n_2}(a)$ divides $\lcm(\alpha_{n_1}(a),\alpha_{n_2}(a))$ by the Chinese remainder theorem.
\end{proof}

Let $n_1$ and $n_2$ be two positive integers such that
$$
\left\{
\begin{array}{ll}
\rad(n_1)|n_2\ \text{and}\ n_2|n_1, & \text{if}\ v_2(n_1)\leq1;\\
2\rad(n_1)|n_2\ \text{and}\ n_2|n_1, & \text{if}\ v_2(n_1)\geq2.
\end{array}
\right.
$$
By definition, we know that $\alpha_{n_1}(a)=\alpha_{n_2}(a)=0$ for every integer $a$ not coprime to $n_1$ and $n_2$. We end this section by determining the exact relationship between $\alpha_{n_1}(a)$ and $\alpha_{n_2}(a)$ for every integer $a$ coprime to $n_1$ and $n_2$.

\begin{thm}\label{thm21}
Let $n_1$ and $n_2$ be two positive integers such that
$$
\left\{
\begin{array}{ll}
\rad(n_1)|n_2\ \text{and}\ n_2|n_1, & \text{if}\ v_2(n_1)\leq1;\\
2\rad(n_1)|n_2\ \text{and}\ n_2|n_1, & \text{if}\ v_2(n_1)\geq2.
\end{array}
\right.
$$
Then, for every integer $a$ coprime to $n_1$ and $n_2$, we have
$$
\alpha_{n_1}(a) = \frac{\alpha_{n_2}(a)}{\gcd\left(\alpha_{n_2}(a),\frac{\gcd(n_1,\mathcal{R}_{n_2}(a))}{n_2}\right)}
$$
\end{thm}

This result is a corollary of the following theorem. 

\begin{thm}\label{thm2}
Let $n_1$ and $n_2$ be two positive integers such that
$$
\left\{
\begin{array}{ll}
\rad(n_1)|n_2\ \text{and}\ n_2|n_1, & \text{if}\ v_2(n_1)\leq1;\\
2\rad(n_1)|n_2\ \text{and}\ n_2|n_1, & \text{if}\ v_2(n_1)\geq2.
\end{array}
\right.
$$
Then, for every integer $a$ coprime to $n_1$ and $n_2$, we have
$$
\ord_{n_1}(a) = \ord_{n_2}(a)\cdot\frac{n_1}{\gcd(n_1,\mathcal{R}_{n_2}(a))}.
$$
\end{thm}

The proof of this theorem is based on the following lemma.

\begin{lem}\label{lem1}
Let $n$ be a positive integer and let $a$ be an integer coprime to $n$. Let $m$ be an integer such that $\rad(m)|\rad{n}$. Then, there exists an integer $u_m$, coprime to $n$ if $m$ is odd, or coprime to $n/2$ if $m$ is even, such that
$$
a^{\ord_n(a)\cdot m} = 1 + u_m\cdot\Rn(a)\cdot m.
$$
\end{lem}

\begin{proof}
We distinguish different cases based upon the parity of $m$. First, we prove the odd case by induction on $m$. If $m=1$, then, by definition of the integer $\Rn(a)$, we have
$$
a^{\ord_n(a)} = 1 + \Rn(a).
$$
Therefore the assertion is true for $m=1$.
\par Now, let $p$ be a prime factor of $m$ and suppose that the assertion is true for the odd number $m/p$, i.e., there exists an integer $u_{m/p}$, coprime to $n$, such that
$$
a^{\ord_{n}(a)\cdot\frac{m}{p}} = 1 + u_{m/p}\cdot\Rn(a)\cdot\frac{m}{p}.
$$
Then, we obtain
$$
a^{\ord_n(a)\cdot m}
\begin{array}[t]{l}
= \displaystyle\left(a^{\ord_n(a)\cdot\frac{m}{p}}\right)^p = \left( 1 + u_{m/p}\cdot\mathcal{R}_n(a)\cdot\frac{m}{p} \right)^p\\[2ex]
= \displaystyle 1 + u_{m/p}\cdot\mathcal{R}_n(a)\cdot m + \sum_{k=2}^{p-1}{\binom{p}{k}\left(u_{m/p}\cdot\mathcal{R}_n(a)\cdot\frac{m}{p}\right)^k} + \left(u_{m/p}\cdot\mathcal{R}_n(a)\cdot\frac{m}{p}\right)^p\\[2ex]
= \begin{array}[t]{l}
\displaystyle 1 + \left( u_{m/p} + \sum_{k=2}^{p-1}{\frac{\binom{p}{k}}{p}\cdot{(u_{m/p})}^{k}\cdot{\Rn(a)}^{k-1}\cdot{\left(\frac{m}{p}\right)}^{k-1}} + \right.\\[2ex]
\displaystyle \left. + {(u_{m/p})}^p\cdot\frac{{\Rn(a)}^{p-1}}{p}\cdot{\left(\frac{m}{p}\right)}^{p-1} \right)\cdot\Rn(a)\cdot m \end{array}\\[2ex]
= \displaystyle 1 + u_m\cdot\Rn(a)\cdot m.
\end{array}
$$
Since $n$ divides $\Rn(a)$ which divides
$$
u_m-u_{m/p} = \sum_{k=2}^{p-1}{\frac{\binom{p}{k}}{p}\cdot{(u_{m/p})}^{k}\cdot{\Rn(a)}^{k-1}\cdot{\left(\frac{m}{p}\right)}^{k-1}} + {(u_{m/p})}^p\cdot\frac{{\Rn(a)}^{p-1}}{p}\cdot{\left(\frac{m}{p}\right)}^{p-1},
$$
it follows that $\gcd(u_m,n)=\gcd(u_{m/p},n)=1$. This completes the proof for the odd case.
\par Suppose now that $n$ and $m$ are even. We proceed by induction on $v_2(m)$. If $v_2(m)=1$, then $m/2$ is odd and by the first part of this proof,
$$
a^{\frac{m}{2}\cdot\ord_n(a)} = 1 + u_{m/2}\cdot\frac{m}{2}\cdot\Rn(a)
$$
where $u_{m/2}$ is coprime to $n$ and hence to $n/2$. Now assume that $v_2(m)>1$ and that
$$
a^{\frac{m}{2}\cdot\ord_n(a)} = 1 + u_{m/2}\cdot\frac{m}{2}\cdot\Rn(a)
$$
with $u_{m/2}$ coprime to $n/2$. Then, we obtain
$$
a^{\ord_n(a)\cdot m}
\begin{array}[t]{l}
= \displaystyle \left(a^{\ord_n(a)\cdot\frac{m}{2}}\right)^2 = \left( 1 + u_{m/2}\cdot\Rn(a)\cdot\frac{m}{2} \right)^2\\[2ex]
= \displaystyle 1 + u_{m/2}\cdot\Rn(a)\cdot m + \left(u_{m/2}\cdot\Rn(a)\cdot\frac{m}{2}\right)^2\\[2ex]
= \displaystyle 1 + \left( u_{m/2} + {(u_{m/2})}^2\cdot\frac{\Rn(a)}{2}\cdot\frac{m}{2} \right)\cdot\Rn(a)m\\[2ex]
= 1 + u_m\cdot\Rn(a)\cdot m.
\end{array}
$$
Since $n/2$ divides $\Rn(a)/2$ which divides $u_m-u_{m/2}$, it follows that $\gcd(u_m,n/2)=\gcd(u_{m/2},n/2)=1$. This completes the proof.
\end{proof}

We are now ready to prove Theorem~\ref{thm2}.

\begin{proof}[Proof of Theorem~\ref{thm2}]
The proof is by induction on the integer $n_1/n_2$. If $n_1=n_2$, then we have
$$
\frac{n_1}{\gcd(n_1,\mathcal{R}_{n_2}(a))} = \frac{n_1}{\gcd(n_1,\mathcal{R}_{n_1}(a))} = \frac{n_1}{n_1} = 1,
$$
since $\mathcal{R}_{n_1}(a)$ is divisible by $n_1$, and thus the statement is true. Let $p$ be a prime factor of $n_1$ and $n_2$ such that $n_2$ divides $n_1/p$ and suppose that
$$
\ord_{n_1/p}(a) = \ord_{n_2}(a)\cdot\frac{n_1/p}{\gcd(n_1/p,\mathcal{R}_{n_2}(a))}.
$$
First, the congruence
$$
a^{\ord_{n_1}(a)} \equiv 1 \pmod{n_1}
$$
implies that
$$
a^{\ord_{n_1}(a)} \equiv 1 \pmod{\frac{n_1}{p}}
$$
and so $\ord_{n_1/p}(a)$ divides $\ord_{n_1}(a)$. We consider two cases.\\[2ex]
\textbf{First Case:} $v_p(n_1)\leq v_p\left(\mathcal{R}_{n_2}(a)\right)$.\\
Since $n_2$ divides $n_1/p$, it follows that $\ord_{n_2}(a)$ divides $\ord_{n_1/p}(a)$. Let $r=\frac{\ord_{n_1/p}(a)}{\ord_{n_2}(a)}$. Hence
$$
\mathcal{R}_{n_1/p}(a)
\begin{array}[t]{l}
= \displaystyle a^{\ord_{n_1/p}(a)}-1 = a^{\ord_{n_2}(a)\cdot r}-1 = \left( a^{\ord_{n_2}(a)}-1 \right)\left( \sum_{k=0}^{r-1}{a^{k\ord_{n_2}(a)}}\right)\\
= \displaystyle\mathcal{R}_{n_2}(a)\left(\sum_{k=0}^{r-1}{a^{k\ord_{n_2}(a)}}\right)
\end{array}
$$
and so $\mathcal{R}_{n_1/p}(a)$ is divisible by $\mathcal{R}_{n_2}(a)$. This leads to
$$
v_{p}(n_1)\leq v_p\left(\mathcal{R}_{n_2}(a)\right)\leq v_p\left(\mathcal{R}_{n_1/p}(a)\right).
$$
Therefore $\mathcal{R}_{n_1/p}(a)$ is divisible by $n_1$ and hence we have
$$
a^{\ord_{n_1/p}(a)} = 1 + \mathcal{R}_{n_1/p}(a) \equiv 1 \pmod{n_1}.
$$
This implies that $\ord_{n_1}(a)=\ord_{n_1/p}(a)$. Moreover, the hypothesis $v_p(n_1)\leq v_p\left(\mathcal{R}_{n_2}(a)\right)$ implies that $\gcd(n_1/p,\mathcal{R}_{n_2}(a))=\gcd(n_1,\mathcal{R}_{n_2}(a))/p$. Finally, we obtain
$$
\ord_{n_1}(a) = \ord_{n_1/p}(a) = \ord_{n_2}(a)\cdot\frac{n_1/p}{\gcd(n_1/p,\mathcal{R}_{n_2}(a))} = \ord_{n_2}(a)\cdot\frac{n_1}{\gcd(n_1,\mathcal{R}_{n_2}(a))}.
$$\ \\
\textbf{Second Case:} $v_p(n_1)>v_p\left(\mathcal{R}_{n_2}(a)\right)$.\\
If $v_2(n_1)\leq1$, then $(n_1/p)/\gcd(n_1/p,\mathcal{R}_{n_2}(a))$ is odd. Otherwise, if $v_2(n_1)\geq2$, then $v_{2}(n_2)\geq2$ and every integer coprime to $n_2/2$ is also coprime to $n_2$. In both cases, $v_2(n_1)\leq 1$ or $v_2(n_1)\geq2$, we know, by Lemma~\ref{lem1}, that there exists an integer $u$, coprime to $n_2$, such that
$$
a^{\ord_{n_1/p}(a)} \begin{array}[t]{l}
= \displaystyle a^{\ord_{n_2}(a)\cdot\frac{n_1/p}{\gcd(n_1/p,\mathcal{R}_{n_2}(a))}} = 1 + u\cdot\mathcal{R}_{n_2}(a)\cdot\frac{n_1/p}{\gcd(n_1/p,\mathcal{R}_{n_2}(a))}\\[2ex]
= \displaystyle 1 + u\cdot\frac{\mathcal{R}_{n_2}(a)}{\gcd(n_1/p,\mathcal{R}_{n_2}(a))}\cdot\frac{n_1}{p}.
\end{array}
$$
As $v_p\left(\mathcal{R}_{n_2}(a)\right)\leq v_p\left(n_1/p\right)$, it follows that $\mathcal{R}_{n_2}(a)/\gcd(n_1/p,\mathcal{R}_{n_2}(a))$ is coprime to $p$, and hence $\ord_{n_1/p}(a)$ is a proper divisor of $\ord_{n_1}(a)$ since
$$
a^{\ord_{n_1/p}(a)} \not\equiv 1 \pmod{n_1}.
$$
Moreover, by Lemma~\ref{lem1} again, there exists an integer $u_p$ such that
$$
a^{\ord_{n_1/p}(a)\cdot p} = 1 + u_p\cdot\mathcal{R}_{n_1/p}(a)\cdot p \equiv 1 \pmod{n_1}.
$$
This leads to
$$
\ord_{n_1}(a) = \ord_{n_1/p}(a)\cdot p = \ord_{n_2}(a)\cdot\frac{n_1}{\gcd(n_1/p,\mathcal{R}_{n_2}(a))} = \ord_{n_2}(a)\cdot\frac{n_1}{\gcd(n_1,\mathcal{R}_{n_2}(a))}.
$$
This completes the proof of Theorem~\ref{thm2}.
\end{proof}

We may view Theorem~\ref{thm2} as a generalization of Theorem~$3.6$ of \cite{Nathanson2000}, where $n_2=p$ is an odd prime number and $n_1=p^k$ for some positive integer $k$. Note that the conclusion of Theorem~\ref{thm2} fails in general in the case where $v_2(n_1)\geq2$ and $n_2=\rad(n_1)$. For instance, for $n_1=24=3\cdot2^3$, $n_2=6=3\cdot2$ and $a=7$, we obtain that $\ord_{n_1}(a)=2$ while $\ord_{n_2}(a)n_1/\gcd(n_1,\mathcal{R}_{n_2}(a))=24/\gcd(24,6)=4$.

We now turn to the proof of the main result of this paper.

\begin{proof}[Proof of Theorem~\ref{thm21}]
From Theorem~\ref{thm2}, we obtain
$$
\alpha_{n_1}(a) \begin{array}[t]{l}
\displaystyle = \ord_{n_1}(a^{n_1}) = \frac{\ord_{n_1}(a)}{\gcd(\ord_{n_1}(a),n_1)} = \frac{\ord_{n_2}(a)\cdot\frac{n_1}{\gcd(n_1,\mathcal{R}_{n_2}(a))}}{\gcd\left(\ord_{n_2}(a)\cdot\frac{n_1}{\gcd(n_1,\mathcal{R}_{n_2}(a))},n_1\right)}\\[2ex]
\displaystyle = \frac{\ord_{n_2}(a)}{\gcd(\ord_{n_2}(a),n_1,\mathcal{R}_{n_2}(a))}.
\end{array}
$$
Thus,
$$
\frac{\alpha_{n_2}(a)}{\alpha_{n_1}(a)}
\begin{array}[t]{l}
= \displaystyle\frac{\frac{\ord_{n_2}(a)}{\gcd(\ord_{n_2}(a),n_2)}}{\frac{\ord_{n_2}(a)}{\gcd(\ord_{n_2}(a),n_1,\mathcal{R}_{n_2}(a))}} = \frac{\gcd(\ord_{n_2}(a),n_1,\mathcal{R}_{n_2}(a))}{\gcd(\ord_{n_2}(a),n_2)}\\[4ex]
= \displaystyle\gcd\left(\frac{\ord_{n_2}(a)}{\gcd(\ord_{n_2}(a),n_2)},\frac{n_2}{\gcd(\ord_{n_2}(a),n_2)}\cdot\frac{\gcd(n_1,\mathcal{R}_{n_2}(a))}{n_2}\right).
\end{array}
$$
Finally, since we have
$$
\gcd\left(\frac{\ord_{n_2}(a)}{\gcd(\ord_{n_2}(a),n_2)},\frac{n_2}{\gcd(\ord_{n_2}(a),n_2)}\right) = \frac{\gcd(\ord_{n_2}(a),n_2)}{\gcd(\ord_{n_2}(a),n_2)} = 1,
$$
it follows that
$$
\frac{\alpha_{n_2}(a)}{\alpha_{n_1}(a)} = \gcd\left(\frac{\ord_{n_2}(a)}{\gcd(\ord_{n_2}(a),n_2)},\frac{\gcd(n_1,\mathcal{R}_{n_2}(a))}{n_2}\right) = \gcd\left(\alpha_{n_2}(a),\frac{\gcd(n_1,\mathcal{R}_{n_2}(a))}{n_2}\right).
$$
\end{proof}

Thus, the determination of $\alpha_n$ is reduced to the case where $n$ is square-free.

\begin{cor}
Let $n$ be a positive integer such that $v_2(n)\leq 1$. Then, for every integer $a$, coprime to $n$, we have
$$
\alpha_n(a) = \frac{\alpha_{\rad(n)}(a)}{\gcd\left(\alpha_{\rad(n)}(a),\frac{\gcd(n,\mathcal{R}_{\rad(n)}(a))}{\rad(n)}\right)}.
$$
\end{cor}

\begin{cor}
Let $n$ be a positive integer such that $v_2(n)\geq 2$. Then, for every integer $a$, coprime to $n$, we have
$$
\alpha_n(a) = \frac{\alpha_{2\rad(n)}(a)}{\gcd\left(\alpha_{2\rad(n)}(a),\frac{\gcd(n,\mathcal{R}_{2\rad(n)}(a))}{2\rad(n)}\right)}.
$$
\end{cor}

\section{The arithmetic function $\beta_n$}

First, we can observe that, by definition of the functions $\alpha_n$ and $\beta_n$, we have
$$
\alpha_n(a) = \beta_n(a) = 0
$$
for every integer $a$ not coprime to $n$ and
$$
\frac{\alpha_n(a)}{\beta_n(a)} \in \left\{1,2\right\}
$$
for every integer $a$ coprime to $n$. There is no general formula known to compute $\alpha_n(a)/\beta_n(a)$ but, however, we get the following proposition.

\begin{prop}\label{prop4}
Let $n_1$ and $n_2$ be two positive integers such that $\rad(n_1)=\rad(n_2)$. Let $a$ be an integer coprime to $n_1$ and $n_2$. If $v_{2}(n_1)\leq 1$, then we have
$$
\frac{\alpha_{n_1}(a)}{\beta_{n_1}(a)} = \frac{\alpha_{n_2}(a)}{\beta_{n_2}(a)}.
$$
If $v_2(n_1)\geq2$, then we have
$$
\alpha_{n_1}(a) = \beta_{n_1}(a).
$$
\end{prop}

\begin{proof}
Let $n_1$ be a positive integer such that $v_2(n_1)\leq 1$ and $a$ be an integer coprime to $n_1$. Let $p$ be an odd prime factor of $n_1$ such that $v_{p}(n_1)\geq2$. We will prove that
$$
\frac{\alpha_{n_1}(a)}{\beta_{n_1}(a)} = \frac{\alpha_{n_1/p}(a)}{\beta_{n_1/p}(a)}.
$$
If $\alpha_{n_1}(a) = 2\beta_{n_1}(a)$, then
$$
{a}^{\beta_{n_1}(a)\cdot n_1} \equiv -1 \pmod{n_1}
$$
and thus
$$
{a}^{\beta_{n_1}(a)\cdot p\cdot\frac{n_1}{p}} \equiv -1 \pmod{\frac{n_1}{p}}.
$$
This implies that $\alpha_{n_1/p}(a) = 2\beta_{n_1/p}(a)$. Conversely, if $\alpha_{n_1/p}(a) = 2\beta_{n_1/p}(a)$, then we have
$$
a^{\beta_{n_1/p}(a)\cdot\frac{n_1}{p}} \equiv -1 \pmod{\frac{n_1}{p}}.
$$
Since $v_p(n_1)\geq2$, it follows that
$$
a^{\beta_{n_1/p}(a)\cdot\frac{n_1}{p}} \equiv -1 \pmod{p}
$$
and thus
$$
{a}^{\beta_{n_1/p}(a)\cdot n_1} + 1 = 1 - {\left( - {a}^{\beta_{n_1/p}(a)\cdot\frac{n_1}{p}} \right)}^p = \left( 1 + {a}^{\beta_{n_1/p}(a)\cdot\frac{n_1}{p}} \right) \sum_{k=0}^{p-1}{{\left( - {a}^{\beta_{n_1/p}(a)\cdot\frac{n_1}{p}} \right)}^k} \equiv 0 \pmod{n_1}.
$$
This implies that $\alpha_{n_1}(a) = 2\beta_{n_1}(a)$. Continuing this process we have
$$
\frac{\alpha_{n_1}(a)}{\beta_{n_1}(a)} = \frac{\alpha_{\rad(n_1)}(a)}{\beta_{\rad(n_1)}(a)}
$$
and since $\rad(n_1)=\rad(n_2)$,
$$
\frac{\alpha_{n_1}(a)}{\beta_{n_1}(a)} = \frac{\alpha_{n_2}(a)}{\beta_{n_2}(a)}.
$$
\par Now, let $n_1$ be a positive integer such that $v_2(n_1)\geq2$, and let $a$ be a non-zero integer. Suppose that we have $\alpha_{n_1}(a)=2\beta_{n_1}(a)$. Since
$$
{a}^{\beta_{n_1}(a)\cdot n_1} \equiv -1 \pmod{n_1}
$$
it follows that
$$
{\left({a}^{\beta_{n_1}(a)\cdot\frac{n_1}{4}}\right)}^{4} \equiv -1 \pmod{4}
$$
in contradiction with
$$
{\left({a}^{\beta_{n_1}(a)\cdot\frac{n_1}{4}}\right)}^{4} \equiv 1 \pmod{4}.
$$
Thus $\alpha_{n_1}(a)=\beta_{n_1}(a)$.
\end{proof}

If $n$ is a prime power, then $\beta_n=\beta_{\rad(n)}$, in analogy with Proposition~\ref{prop3} for $\alpha_n$.

\begin{prop}
Let $p$ be a prime number and let $a$ be an integer. Then we have
$$
\beta_{p^k}(a) = \beta_p(a)
$$
for every positive integer $k$.
\end{prop}

\begin{proof}
This result is trivial for every integer $a$ not coprime to $p$. Suppose now that $a$ is coprime to $p$. For $p=2$, then, by Proposition~\ref{prop4}, we have
$$
\beta_{2^k}(a) = \alpha_{2^k}(a) = 1
$$
for every positive integer $k$. For an odd prime number $p\geq3$, Proposition~\ref{prop4} and Proposition~\ref{prop3} lead to
$$
\beta_{p^{k}}(a) = \frac{\alpha_{p^k}(a)}{\alpha_p(a)}\cdot\beta_p(a) = \beta_p(a)
$$
for every positive integer $k$. This completes the proof.
\end{proof}

Let $n_1$ and $n_2$ be two positive integers such that
$$
\left\{
\begin{array}{ll}
\rad(n_1)|n_2\ \text{and}\ n_2|n_1, & \text{if}\ v_2(n_1)\leq1;\\
2\rad(n_1)|n_2\ \text{and}\ n_2|n_1, & \text{if}\ v_2(n_1)\geq2.
\end{array}
\right.
$$
It immediately follows that $\beta_{n_1}(a)=\beta_{n_2}(a)=0$ for every integer $a$ not coprime to $n_1$ and $n_2$. Finally, we determine the relationship between $\beta_{n_1}(a)$ and $\beta_{n_2}(a)$ for every integer $a$ coprime to $n_1$ and $n_2$.

\begin{thm}\label{thm31}
Let $n_1$ and $n_2$ be two positive integers such that
$$
\left\{
\begin{array}{ll}
\rad(n_1)|n_2\ \text{and}\ n_2|n_1, & \text{if}\ v_2(n_1)\leq1;\\
2\rad(n_1)|n_2\ \text{and}\ n_2|n_1, & \text{if}\ v_2(n_1)\geq2.
\end{array}
\right.
$$
Let $a$ be an integer coprime to $n_1$ and $n_2$. Then, we have
$$
\beta_{n_1}(a) = \frac{\beta_{n_2}(a)}{\gcd\left(\beta_{n_2}(a),\frac{\gcd(n_1,\mathcal{R}_{n_2}(a))}{n_2}\right)}.
$$
\end{thm}

\begin{proof}
If $v_2(n_1)\leq 1$, then Theorem~\ref{thm21} and Proposition~\ref{prop4} lead to
$$
\frac{\beta_{n_2}(a)}{\beta_{n_1}(a)} = \frac{\alpha_{n_2}(a)}{\alpha_{n_1}(a)} = \gcd\left(\alpha_{n_2}(a),\frac{\gcd(n_1,\mathcal{R}_{n_2}(a))}{n_2}\right).
$$
Since $v_2(n_2)=v_2(n_1)\leq1$, it follows that $\gcd(n_1,\mathcal{R}_{n_2}(a))/n_2$ is odd and hence, we have
$$
\frac{\beta_{n_2}(a)}{\beta_{n_1}(a)} = \gcd\left(\alpha_{n_2}(a),\frac{\gcd(n_1,\mathcal{R}_{n_2}(a))}{n_2}\right) = \gcd\left(\beta_{n_2}(a),\frac{\gcd(n_1,\mathcal{R}_{n_2}(a))}{n_2}\right).
$$
If $v_{2}(n_1)\geq2$, then $\beta_{n_1}(a)=\alpha_{n_1}(a)$ and $\beta_{n_2}(a)=\alpha_{n_2}(a)$ by Proposition~\ref{prop4} and the result follows from Theorem~\ref{thm21}.
\end{proof}

Thus, as for $\alpha_n$, the determination of $\beta_n$ is reduced to the case where $n$ is square-free.

\begin{cor}
Let $n$ be a positive integer such that $v_2(n)\leq 1$. Then, for every integer $a$, coprime to $n$, we have
$$
\beta_n(a) = \frac{\beta_{\rad(n)}(a)}{\gcd\left(\beta_{\rad(n)}(a),\frac{\gcd(n,\mathcal{R}_{\rad(n)}(a))}{\rad(n)}\right)}.
$$
\end{cor}

\begin{cor}
Let $n$ be a positive integer such that $v_2(n)\geq 2$. Then, for every integer $a$, coprime to $n$, we have
$$
\beta_n(a) = \frac{\beta_{2\rad(n)}(a)}{\gcd\left(\beta_{2\rad(n)}(a),\frac{\gcd(n,\mathcal{R}_{2\rad(n)}(a))}{2\rad(n)}\right)}.
$$
\end{cor}

\section{Acknowledgments}

The author would like to thank Shalom Eliahou for his help in preparing this paper. He also thanks the anonymous referee for its useful remarks.

\nocite{*}
\bibliographystyle{plain}
\bibliography{biblio}

\end{document}